\documentclass[journal,onecolumn,draftcls]{IEEEtran}

	\usepackage{url}
	\usepackage{ifpdf}
	\ifpdf
	\usepackage[pdftex]{graphicx}
	\else
	\usepackage[dvips]{graphicx}
	\fi

	\DeclareGraphicsExtensions{.eps}
	\usepackage{epstopdf}
	\usepackage{tikz}

	\usetikzlibrary{shapes.geometric,arrows}
	\tikzstyle{decision} = [diamond, draw, fill=white!20,
	text width=4.5em, text badly centered, node distance=2cm, inner sep=0pt]
	\tikzstyle{block} = [rectangle, draw, fill=white!20,
	text width=15em, text centered, rounded corners, node distance=3cm, minimum height=4em]
	\tikzstyle{smallblock} = [rectangle, draw, fill=white!20,
	text width=5em, text centered, rounded corners, node distance=3cm, minimum height=4em]
	\tikzstyle{line} = [draw, -latex']
	\tikzstyle{cloud} = [draw, ellipse,fill=white!20,  node distance=2cm,
	minimum height=2em]

	\usepackage{subcaption}
	\usepackage{caption}
	\usepackage{epsfig}
	\usepackage{amsmath}
	\usepackage{url}
	\usepackage{color}
	\usepackage{algorithm}
	\usepackage{algpseudocode}
	\usepackage{mathcomp}
	\usepackage{mathtools}
	\usepackage{amssymb}
	\usepackage{amsmath}
	\usepackage{amsthm}
	\usepackage{multicol}
	\usepackage{multirow}

	\newtheorem{theorem}{Theorem}
	\newtheorem*{rep@theorem}{\rep@title}
	\newcommand{\newreptheorem}[2]{%
		\newenvironment{rep#1}[1]{%
			\def\rep@title{#2 \ref{##1}}%
			\begin{rep@theorem}}%
			{\end{rep@theorem}}}
	\makeatother
	\newreptheorem{theorem}{Theorem}

	\newtheorem{definition}{Definition}
	
	\newtheorem{remark}{Remark}
	\newtheorem{proposition}{Proposition}
	\usepackage{mathrsfs}
	\usepackage{graphicx}
	\usepackage{hhline}
	\usepackage{bigstrut}
	\usepackage{cases}
	\usepackage{verbatim}
	\usepackage{bm}

	\DeclareMathOperator{\E}{\mathbb{E}}
	\DeclareMathOperator{\Var}{\text{Var}}
	\DeclareMathOperator{\Cov}{\text{Cov}}

\begin{document}
\title{A Distributed Online Pricing Strategy for Demand Response Programs}

\author{Pan~Li,~\IEEEmembership{Student~Member,~IEEE,}~Hao Wang,~\IEEEmembership{Member,~IEEE}~and~Baosen Zhang,~\IEEEmembership{Member,~IEEE}
\thanks{This work was partially supported by NSF grant CNS-1544160 and the University of Washington Clean Energy Institute.}
\thanks{The authors are with the Department of Electrical Engineering, University of Washington, Seattle, WA 98195, USA (e-mail: \{pli69,hwang16,zhangbao\}@uw.edu).}}

\markboth{IEEE Transactions on Smart Grid}%
{Li\MakeLowercase{\textit{et al.}}: TBD}

\maketitle
	\begin{abstract}
We study a demand response problem from utility (also referred to as operator)'s perspective with realistic settings, in which the utility faces uncertainty and limited communication. Specifically, the utility does not know the cost function of consumers and cannot have multiple rounds of information exchange with consumers. We formulate an optimization problem for the utility to minimize its operational cost considering time-varying demand response targets and responses of consumers. We develop a joint online learning and pricing algorithm. In each time slot, the utility sends out a price signal to all consumers and estimates the cost functions of consumers based on their noisy responses. We measure the performance of our algorithm using regret analysis and show that our online algorithm achieves logarithmic regret with respect to the operating horizon. In addition, our algorithm employs linear regression to estimate the aggregate response of consumers, making it easy to implement in practice. Simulation experiments validate the theoretic results and show that the performance gap between our algorithm and the offline optimality decays quickly.

	\end{abstract}


\IEEEpeerreviewmaketitle

	\section{Introduction}\label{intro}
	This paper considers an important problem in demand response programs. Demand response (DR) is a popular mechanism used by utilities to modify the electricity consumption of consumers and has been extensively studied by the community in recent years~(see, e.g.~\cite{PalenskyEtAl2011,AlbadiEtAl2008} and the references within).  When a utility (also referred to as operator) engages its customer in demand response, it typically sends out a signal--for example, price~\cite{NguyenEtAl2011}, incentives~\cite{ZhongEtAl2013}, or even text messages~\cite{ENERNOC2017}-- to elicit changes in electricity consumption. 

In practice, the utility faces two key challenges in these type of programs.\footnote{We do not consider direct load control for demand response in this paper~(e.g., see~\cite{ChenEtAl2014}).} The first is uncertainty, where utilities often do not know the exact response of the customers to the signals. Secondly, utilities often have specific demand response targets at a particular time. For example, if utility wishes to reduce the demand of a residential area by $100$ KW, it is suboptimal to receive either $50$ KW or $200$ KW from a DR call. Thus DR programs are faced with the fundamental problem of optimizing the responses under uncertainty about its customers. In this paper, we address this problem by designing optimal algorithms that can optimize and learn simultaneously.

The problem of private customer information is not new in demand response. A standard approach to address this challenge is to adopt a ``negotiation'' process between the customers and the utility. For example, at time $t$, suppose the utility wishes to elicit some change from the customers but do not know their individual cost functions. By repeatedly probing the customer with different prices, the optimal price--i.e., same price as if the utility had full information--can be agreed upon as long as the cost functions satisfy some technical constraints~\cite{LiEtAl2011}. However, in practice, this strategy may be difficult to implement in some settings. Firstly, some systems still lack real-time two-way communication between customers and utilities, which makes multiple rounds of information exchange in a short time difficult. Secondly, many customers do not actively engage in the process and tend to change their consumption just based on the first price they observe. For example, many companies use text-based notifications~\cite{ENERNOC2017,Ohmconnect2017}, and it is unrealistic to assume customers would manually negotiate with the companies. In this paper, we move away from this setting by assuming that at any time $t$, a customer changes its consumption after it receives a price and does not undergo a negotiation process. Therefore any information the utility learns from the price and customers' responses at time $t$ will only be useful the next time demand response is called. Thus some loss relative to the case where the utility has full information, or \emph{regret}, is inevitable. Our contribution is in describing an algorithm that minimizes the regret.

Specifically, we consider the following optimization problem. An utility has a \emph{time-varying} demand response target at each time $t$ and manages $N$ users. Each of the users has a quadratic cost function that is not known to the utility. To achieve its desired target, the utility sends out a price signal to the users. Interpreting this price as payment, the users solve an optimization problem to determine their changes in demand. The utility observes the response of the customers. The two key difficulties of the problem are: 1) the utility \emph{commits} to the price at the beginning of a time slot and can only use its observations for the next time; and 2) the demand response target changes throughout time. Therefore the price must be determined in an online fashion, optimized using historical prices and responses.

To measure the performance of an online problem, we adopt the standard metric of \emph{regret}, which measures the difference between the performance of the online scheme and the performance of an optimal offline scheme that has the full information. Normally regret is reported as the sum of the differences between the online and offline versions as a function of the total number of time steps $T$. A regret \emph{sublinear in $T$} is consider to be a ``good'' regret, since a trivial online algorithm has regret that is proportional to $T$. In a seminal work, the authors in \cite{LaiEtAl1985} showed that under broad conditions, no online algorithm can achieve a better regret than $\log(T)$. In this paper, we present an algorithm that meets this lower bound. Throughout the paper, we assume that the utility has a quadratic penalty if it receives a different demand response amount than its target and each user has a convex quadratic cost function. We do not explicitly model power flow, and view the aggregate demand response received by the utility as the sum of the response of each user. We believe losses can be include into the framework by accounting for the topology of the power network.

\subsection{Literature Review} \label{sec:lit}
Demand response was first proposed as a load control mechanism and gained popularity in both academia and industry in the past decade because of the advancement of smartgrid~\cite{AlbadiEtAl2008,PalenskyEtAl2011,Qdr2006}. The setting of where utility does not have (or only has partial) customer information was considered by authors in \cite{SaadEtAl2012,LiEtAl2011,VardakasEtAl2015}. Most of these studies assumed that the utility and users conduct a negotiation process before committing to a decision.  This process is typically modeled as implementing a distributed optimization algorithm~\cite{TsuiEtAl2012,JiangEtAl2011,SamadiEtAl2013}.

Recently, a series of papers investigated the setting where the utility must commit to its decision under uncertainty about its users. The authors in~\cite{ZhangEtAl2016} considered using reinforcement learning to learn the home appliances of a single household, and the related work of designing an automated system in~\cite{ONeillEtAl2010} include energy scheduling under uncertainty in price. The authors in~\cite{XuEtAl2016} considered a data-driven approach to learn customer behavior using historical consumption profile. The authors in~\cite{BaltaogluEtAl2016} considered a two-stage market where the demand is modeled by a Markov jump process. Our work differs from these in that we focus on the problem from the utility's point of view, focus on the uncertainty in the response of the customers, and do not require extensive data on historical consumption profiles. The work in~\cite{KhezeliEtAl2016} is closest to our model, although under a different setup where the authors consider risk in learning and not directly focused on cost of the utility.


\subsection{Contribution of this paper}
We state the main contributions of this paper as the following:
\begin{itemize}
	\item We formulate a framework that models the practical situation where utility calls on demand response under time-varying demand response requirements without full information about its customers.
	\item We propose an online algorithm that achieves the optimal regret of $\log(T)$, where $T$ is the total number of time periods. This algorithm leverages the fact that the utility only needs to know the \emph{aggregate} response of it's customers, making the learning process scalable. We validate this property via both theoretic analysis and simulations.
	\item The algorithm employs standard least square estimation techniques to find the correct price, makes it simple to apply in practice. The standard estimation technique exhibits fast convergence under various scenarios. In addition, our analysis clearly illustrates the impact of uncertainty about customers cost functions on the regret. 
\end{itemize}

The remainder of this paper is organized as follows. We formulate the demand response optimization problem from the operator's perspective and state the main results in Section \ref{formulation}. We then analyze the offline optimal solution in Section \ref{offline}. In Section \ref{online}, we design the online learning algorithm using linear regression to estimate the unknown parameters of users and show that the algorithms achieve performance of $\log (T)$-regret. Numerical results are presented in Section \ref{simulation}. Finally, we conclude our work in Section \ref{conclusion}.

	\section{Model and Problem Formulation}\label{formulation}
\begin{figure}[!htbp] 
	\centering
	\includegraphics[width = 0.35\columnwidth]{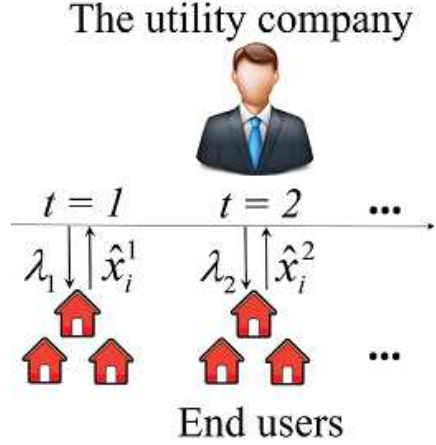}
	\caption{Interactions between the utility and users: $\lambda_t$ is the price at time $t$ and $\hat{x}_i^t$ is the change in demand of user $i$ at time $t$.}
	\label{fig-system}
	\end{figure}

	We consider a demand response system shown in Fig. \ref{fig-system}.
An utility coordinates a set of users $i=1,...,N$ to provide demand response service to the power grid over a time horizon $t = 1,...,T$. At time $t$, the utility broadcasts a price $\lambda_t$ to all users, and each user $i$ responds to the price and changes its consumption by some amount.
We assume that there are no iterative communications between users and the utility in a single time period. That is, a user responds at the time he/she receives the price signal without any negotiations. In the following, we will present the models for both utility and users and formulate the demand response optimization problem and its offline optimal solution.

\subsection{Utility Cost Model}

The utility aggregates power consumption changes from all the users to provide demand response service to the power grid \cite{LiEtAl2011}. The utility seeks to minimize its cost or maximize its profits by calling a certain amount of demand responses from users. We assume that the power grid specifies a normalized amount of demand reduction denoted as $d_t$ at the beginning of time slot $t$. The utility determines its capacity of demand response $Y$ and provides a response of $Y d_t$. We assume that the utility receives revenue proportional to a price $\alpha $. The capacity $Y$ captures size of the utility, that is, the number of customers an utility can call to provide demand response.
However, if the actual aggregated demand reduction deviates from the committed reduction $Y d_t$, the utility will suffer a penalty as its cost. Let ${x}_i^t$ denote the response of user $i$ at time $t$. Then the utility's profit in time slot $t$ is modeled as
\begin{equation}
	\alpha Y - \frac{1}{2} \left(\sum_{i=1}^{N} {x}_i^t - Y d_t \right)^2.
\end{equation}

The first term $\alpha Y$ is the revenue from providing demand reduction at capacity $Y$. where $\alpha$ is the unit price of a standardized size of demand reduction (denoted by $d_t$) and $Y$ is a scaling factor showing how many standardized demand reduction the utility provides. The second term is the penalty of deviation between actual aggregated demand reduction $\sum_{i=1}^{N} {x}_i^t$ from the committed reduction $Y d_t$.

\subsection{User Model}
We model users' response to the broadcasted price as solving an optimization problem to balance their reward from responding with respect to their cost functions~\cite{AlbadiEtAl2008,LiEtAl2011}. Here we consider the following quadratic cost function $u_i(\cdot)$ for user $i$:
\begin{equation}\label{consumer_model}
	u_i(x_i^t) = \frac{1}{2} \beta_i (x_i^{t})^2 + \alpha_i x_i^{t},
\end{equation}
where $x_i^t$ is user $i$'s response in consumption at time $t$.

Note that both $\beta_i$ and $\alpha_i$ are unknown parameters to the utility.
The price signal serves as a reward to users for their consumption reduction. Specifically, for those users who adjust their consumption by $x_i^{t}$, they will receive payment $\lambda'_t x_i^{t}$ from the utility in time slot $t$, where $\lambda'_t $ is the price in unit consumption. Since users aim to minimize their net cost from reducing consumption, we have the cost minimization problem for each user $i$:
\begin{equation}\label{user_optimization}
\begin{aligned}
& \min_{x_i^t} && u_i(x_i^t) - \lambda'_t x_i^t .
\end{aligned}
\end{equation}

Note that this optimization problem is operationally the same as $\min_{x_i^t}\ \frac{u_i(x_i^t)}{N} - \lambda_t x_i^t $, where $\lambda_t = \frac{\lambda_t'}{N}$. This is in accordance with the optimization that is introduced later in Section \ref{sec:prob_form}, when the costs are divided by $N$. To be consistent, we refer to $\lambda_t$ as the price signal that stimulates consumption. 

In practice, there is always some noise associated with the response of users. For example, users~(or more likely their home management systems) can relax the temperature tolerance of their homes, but the actual power consumption reduction would be somewhat random~\cite{AlbertEtAl2015a}.
Therefore, we assume that the actual response of the users is a noisy version of the optimal solution to \eqref{user_optimization}. 
Then user $i$'s consumption change is given by
\begin{equation} \label{optimal2_noise0}
\hat{x}_{i}^{t} = x_{i}^{t} + \epsilon_i^t,
\end{equation}
where $\epsilon_i^t$ is the noise in observation and is independent and identically distributed. We assume the noise is Gaussian in this paper, although the results extend to other types of noises.

	\subsection{Problem Formulation}\label{sec:prob_form}
The operator aims to minimize the following expected average cost over the time horizon:
\begin{equation}\label{case2_formulation}
\begin{aligned}
& \min_{Y}~\min_{\bm{x}}  \sum_{t=1}^{T}  \sum_{i=1}^{N} \frac{1}{N} \E \left[ \frac{1}{2} \beta_i (x_i^{t}+\epsilon_i^t)^2 + \alpha_i (x_i^{t}+\epsilon_i^t) \right] \\
& + \sum_{t=1}^{T} \frac{1}{2N} \E\left[\left(\sum_{i=1}^{N} (x_i^t+\epsilon_i^t) - Y d_t \right)^2\right] - \frac{\alpha Y T}{N} \\
= & \min_{Y}~\min_{\bm{x}} \E  \sum_{t=1}^{T} C(\bm{\hat{x}}_t,Y) - \frac{\alpha Y T}{N}
\end{aligned}
\end{equation}
where the expectation is taken over the noise $\epsilon_i^t$ in user's response. The cost function $C(\bm{\hat{x}}_t,Y)$ involves two terms: the quadratic cost of users' response in consumption, i.e.,  $\bm{\hat{x}}_t$, and the cost of mismatch between the sum of user consumption and the target consumption. Here, $\bm{\hat{x}}_t$ is a vector storing each user's consumption $\hat{x}_i^t$ at time $t$. The last term $\frac{\alpha Y T}{N}$ is the revenue of the utility. 
The utility determines the optimal capacity of demand response $Y$ at the beginning of the time horizon, and the optimal responses $x_i^t$'s are determined in every time period. Since $\lambda^t$ uniquely determines the responses~(up to the random noises), we sometimes write them as $\hat{x}_i^t(\lambda^t)$.

The optimization problem \eqref{case2_formulation} is not difficult to solve if $\beta_i$ and $\alpha_i$ are known to the utility. It suffices to solve a convex quadratic programming in a centralized manner. Alternatively, if the users can iteratively communicate with the utility in a time step, a decentralized algorithm can also be deployed. However, utilities generally do not have information about the cost coefficients of the users~\cite{YangEtAl2014, YangEtAl2015}, and often do not engage in multi-round communication with users~\cite{KhezeliEtAl2016}. In the next section, we present the main results of this paper, which describes how \eqref{case2_formulation} can be solved when the utility needs to learn and optimize demand response of users in an online fashion.




\subsection{Online Algorithm and Regret}
%

When the utility does not know the parameters of the cost function and wants to solve the optimization problem in \eqref{case2_formulation}, it needs to solve two problems at the same time: 1) learn the cost function, and 2) design proper price signals to obtain desired aggregated responses from end users. These two problems are solved in an \emph{online fashion}, where the utility determines current price ${\lambda}_t$ based on partial information, i.e., the past sequence of prices $\{{\lambda}_1, {\lambda_2}, ..., {\lambda}_{t-1} \}$, and the past responses from users $\{\hat{x}_i^s(\lambda_s), 1\leq i \leq N, 1 \leq s \leq t-1 \}$ . The online strategy is a policy that maps the history observation $\mathcal{H}_t = \{{\lambda}_s, \hat{x}_i^s(\lambda_s), 1\leq i \leq N, 1 \leq s \leq t-1 \}$ to ${\lambda}_t$, denoted by ${\lambda}_t = f(\mathcal{H}_t)$. The function $f(\cdot)$ is called a policy, or more specifically, an online learning algorithm that uses information from the past to make sequential decisions without knowing the future information. The main result is that there is an online policy that will eventually \emph{find the correct price}, with a small regret.

We adopt the notion of regret to evaluate the performance of an online pricing strategy $\lambda_t$. Regret is widely used in literature to evaluate an online decision process, for example in online convex problems \cite{Shaley2012} and in multi-arm bandit problems~\cite{Bubeck2012}. In this paper, it is defined as the expected difference between the costs obtained by the online strategy with unknown $\alpha_i$ and $\beta_i$, and the optimal solution of when all parameters are known in advance. More specifically, regret $R$ over horizon $1, 2, ..., T$ is defined as:
\begin{equation}\label{regret}
\begin{aligned}
R
& = \E \{\sum_{t = 1}^T C(\bm {\hat{x}}_t(\lambda_{t}), Y)\}  - \E \{\sum_{t = 1}^T C(\bm{ \hat{x}}_t(\lambda^*_{t}), Y)\},
\end{aligned}
\end{equation}
where $Y$ is the capacity predetermined beforehand for the demand response program. The expectation is taken over the randomness from $\bm {\hat{x}}_t(\lambda_{t})$, when users' consumption is noisy.

 In addition, $\lambda^*_{t}$ is the optimal reward signal if given full information (values of $\alpha_i$, $\beta_i$), and $\lambda_{t}$ is the online pricing signal obtained from online strategy. We stack all users' response $x_i^t$ as a vector $\bm {\hat{x}}_t(\lambda_{t})$ at time $t$ stimulated with price $\lambda_{t}$.
The regret can be decomposed as across each time period as:
\begin{equation}\label{regrett}
\begin{aligned}
R
& = \sum_{t = 1}^T \{ \E \{ C(\bm {\hat{x}}_t(\lambda_{t}), Y)\}  - \E \{ C(\bm{ \hat{x}}_t(\lambda^*_{t}),Y)\} \} \\
& = \sum_t R_t . \\
\end{aligned}
\end{equation}


\subsection{Summary of main results} \label{mainresults}
We now state the main theorem in this paper, that quantifies the performance of the online strategy. The details of the learning strategy are illustrated in Section \ref{online}.


The main theorem of this paper is:
\begin{theorem}\label{maintheorem}
	There is an algorithm where the gap $R_t$ between online strategy and optimal offline solution decays as $1/t$ and the algorithm achieves $\log T$-regret.
\end{theorem}

\begin{figure}[!ht]
	\centering
	\includegraphics[width = 0.9\columnwidth]{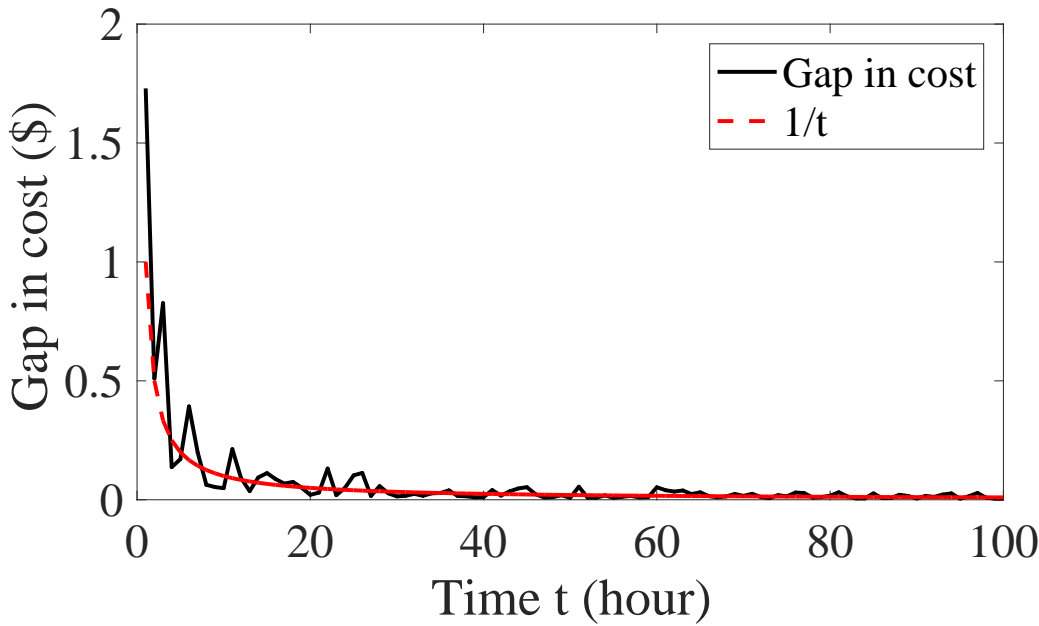}
	\caption{Gap between the online cost and optimal cost over time, accumulated over 100 users for 100 time slots. Detailed parameters in the model are illustrated in Section \ref{simulation}.}
	\label{eg1}
\end{figure}

Theorem \ref{maintheorem} states that the gap $R_t$ is diminishing with time $t$ and has a rate $\frac{1}{t}$, as shown in Fig. \ref{eg1}. This gives an overall regret in the order of $\log(T)$ over at total of $T$ time periods. This is actually the best one can hope for, because of a celebrate result in online optimization that states no online algorithm can do better~\cite{LaiEtAl1985}.

In Section \ref{online}, we aim to find an algorithm such that this bound of the regret is achieved. The existence of such an algorithm validates the statement in Theorem \ref{maintheorem}, and the proof on the convergence of the algorithm is left to the appendices for interested readers.



	\section{Offline optimal solution}\label{offline}
In this section, we characterize the offline optimal solution, which provides insights for designing the online strategy in Section \ref{online}.

Suppose that the utility knows the value of $\alpha_i$ and $\beta_i$ of all users. To solve \eqref{case2_formulation}, we start by solving the inner minimization problem, assuming that the optimal capacity $Y^*$ has been determined beforehand. Then we solve the following problem to derive the optimal demand response of all the users $i = 1, 2, ..., N$ over time slots $t = 1, 2, .., T$, with $Y^*$ as a parameter:

\begin{equation}\label{simplified0}
\begin{aligned}
\min_{\bm{x}} & \sum_{t=1}^{T}  \sum_{i=1}^{N} \frac{1}{N} \E \left( \frac{1}{2} \beta_i (x_i^{t}+\epsilon_i^t)^2 + \alpha_i (x_i^{t}+\epsilon_i^t) \right) \\
&+ \sum_{t=1}^{T}  \frac{1}{2N} \E\left[(\sum_{i=1}^{N} x_i^t +\epsilon_i^t - Y^* d_t)^2\right] \\ 
\end{aligned}
\end{equation}


Since $Y^*$ is a predefined scaling factor, each time period is decoupled.  The optimal scheduling for each user $i$ at each time $t$ is computed as the solution from \eqref{simplified0}:
\begin{align} \label{optimal0}
& x_{i}^{t,*} = \frac{N \frac{Y^* d_t + \sum_{i=1}^N \frac{\alpha_i}{\beta_i} }{N+ \sum_{i=1}^N \frac{N}{\beta_i}} - \alpha_i}{\beta_i},
\end{align}
where we used the fact that the noise is zero mean and independent to everything~\cite{Chow76}.

In the online setting, the utility can only influence the users through a price. Therefore here we look at the price signals that would realize the desired changes in \eqref{optimal0}.
To design this price, we introduce an ancillary variable $Q_t$, which is the sum of responses from users into problem in \eqref{simplified0}:
\begin{equation}\label{simplified}
\begin{aligned}
\min_{\bm{Q}, \bm{x} } & \sum_{t=1}^{T}  \sum_{i=1}^N \frac{1}{N} \E \left( \frac{1}{2} \beta_i (x_i^{t}+\epsilon_i^t)^2 + \alpha_i (x_i^{t}+\epsilon_i^t) \right) \\
&+ \sum_{t=1}^{T} \frac{1}{2N} \E\left[\left(Q_t+\sum_{i=1}^N \epsilon_i^t - Y^* d_t\right)^2\right] \\ 
\text{s.t.} \ & Q_t- \sum_{i=1}^{N} x_i^t = 0.
\end{aligned}
\end{equation}

The problem in \eqref{simplified} is a constrained convex problem \cite{boyd2004convex} that solves exactly the same problem as presented in \eqref{simplified0}. We introduce a dual variable $\lambda_t$ for the equality constraint
$Q_t- \sum_{i\in\mathcal{N}} x_i^t = 0$ and use KKT conditions \cite{boyd2004convex} to solve the optimal solution in a closed form as:
\begin{equation}\label{optimal2}
	x_{i}^{t,*} = x_i^t(\lambda_t^*) = \frac{N \lambda_t^* - \alpha_i}{\beta_i},
\end{equation}
and
\begin{equation}\label{optimallambda}
 \lambda_{t}^{*} = \frac{Y^* d_t +\sum_{i=1}^N \frac{\alpha_i}{\beta_i} }{N+ \sum_{i=1}^N \frac{N}{\beta_i}},
\end{equation}
where $\lambda_{t}^{*}$ is the optimal value of the dual variable associated with the equality constraint in \eqref{simplified} and can be interpreted as optimal price incentive that the utility needs to incentivize users to achieve the optimal consumption change as in centralized control in \eqref{simplified0}. Note that since each cost term in \eqref{simplified0} is divided by a factor of $N$, $\lambda_t$ and $\lambda_t'$ (defined in \eqref{user_optimization}) differs with a scaling factor of $\frac{1}{N}$ as well. 



After solving the optimal demand change of users $\bm{x}^{*} = \{ x_{i}^{t,*},~\forall i,~\forall t \}$, we substitute the optimal solution obtained in \eqref{optimal0} as a function of $Y$ into \eqref{simplified0} and solve the optimal solution for $Y^*$ as
\begin{equation} \label{Ystar}
Y^* = \frac{T\alpha(1+\sum_i \frac{1}{\beta_i})^2 - \sum_t (\sum_i \frac{\alpha_i}{\beta_i})(1+ \sum_i \frac{1}{\beta_i})d_t}{\sum_t d_t^2 (1+\sum_i \frac{1}{\beta_i})},
\end{equation}
which depicts the optimal ``capacity'' in demand reduction during demand response programs.




%


	\section{Online strategy and regret analysis}\label{online}
	

Now we proceed to design an optimal online pricing strategy that achieves the regret as mentioned in Theorem \ref{maintheorem}, given that the parameters are not unknown to the utility. There are two questions: 1) does the utility need to estimate each individual $\alpha_i$ and $\beta_i$ to design such a price, and 2) does the noisy observation from the users' consumption data affect the estimating procedure.


First, given full information, the optimal price signal presented in \eqref{optimallambda} only involves the terms $\sum_i \frac{1}{\beta_i}$ and $\sum_i \frac{\alpha_i}{\beta_i}$, given that $Y^*$ is set to some predefined value. This suggests that if $\alpha_i$ and $\beta_i$ are not known to the utility, only those two aggregated terms $\sum_i \frac{1}{\beta_i}$ and $\sum_i \frac{\alpha_i}{\beta_i}$ need to be computed. The dimension of the parameter space reduces from $2N$ to $2$, which reduces the information needed to compute the online price. More specifically, it is sufficient for the utility to observe the aggregated consumption:
\begin{equation}\label{agg}
	\sum_i \hat{x}_i^t = \sum_i \frac{N}{\beta_i} \lambda_t - \sum_i\frac{\alpha_i}{\beta_i} + \sum_i\epsilon_i^t
\end{equation}
where we restate that $\lambda_t$ is the designed price signal at time $t$ that is based on past observation. Since the noise for each user is independent and follows standard gaussian distribution, the aggregated noise $\sum_i\epsilon_i^t $ follows a gaussian distribution with variance equals to $N$.

Second, note that the observation of the aggregated consumption in \eqref{agg} follows a gaussian distribution and the expected value of the observation is linear in the parameter $\sum_i \frac{N}{\beta_i}$ and $\sum_i\frac{\alpha_i}{\beta_i}$ . Since least square estimator is the Best Linear Unbiased Estimator (BLUE) \cite{BLUE}, we therefore adopt the least square approach to estimate $\sum_i \frac{N}{\beta_i}$ and $\sum_i\frac{\alpha_i}{\beta_i}$ using a linear regression model.

With the clarification of the two questions, we elaborate the formulation of the regret introduced in Section \ref{offline} and propose an iterative linear regression estimator in Algo. \ref{lin1} to design online price strategy. In the following subsections, we first give more details on the formulation of the regret. Then we introduce the proposed online learning strategy for pricing $\lambda_t$ based on historical observation. We also state the performance of the online strategy and leave the detailed analysis on the order of the regret (Theorem \ref{maintheorem2}) to the appendices.





\subsection{Regret formulation}

As stated in Section \ref{formulation}, we evaluate an online pricing strategy by regret. We let $\hat{x}_i^{t*} = {x}_i^{t*} + \epsilon_i^t$ and $\hat{Q}_t^* = \sum_i \hat{x}_i^{t*}$ denote the optimal response with noise and optimal aggregated response, respectively. Since each time $t$ is decoupled and for simplicity reasons, we show the computation for one-step regret (gap $R_t$) as the following:

\begin{equation} \label{regretcase2}
\begin{aligned}
R_t =& \E \{ C(\bm {\hat{x}}_t(\lambda_{t}), Y^*)\}  - \E \{ C(\bm{ \hat{x}}_t(\lambda^*_{t}), Y^*)\} \\
=& \E \left( \frac{\sum_{i=1}^N \beta_i (\hat{x}_i^t)^2 + 2 \alpha_i \hat{x}_i^t}{2N} + \frac{1}{2N} \left( \hat{Q}_t - Y^* d_t \right)^2 \right) \\
& - \left( \frac{\sum_{i=1}^N \beta_i ( \hat{x}_i^{t*})^2 + 2\alpha_i \hat{x}_i^{t*}}{2N} + \frac{1}{2N} \left( \hat{Q}_t^* - Y^* d_t \right)^2 \right) \\
=& \frac{N}{2} \left(\sum_{i=1}^N  \frac{1}{\beta_i} + (\sum_{i=1}^N \frac{1}{\beta_i} )^{2} \right)
\left[ \E \lambda_{t}^{2} - (\E \lambda_{t})^{2} + (\E \lambda_{t}- \lambda_{t}^{*})^{2}  \right] \\
&+ \sum_{i=1}^N \frac{\alpha_i}{\beta_i} \left( \sum_{i=1}^N \frac{1}{\beta_i} -1 \right) (\E \lambda_{t} - \lambda_{t}^{*} ) \\
=& C_1 \left[ \E \lambda_{t}^{2} - (\E \lambda_{t})^{2} + (\E \lambda_{t}- \lambda_{t}^{*})^{2}  \right] + C_2 (\E \lambda_{t} - \lambda_{t}^{*} ), \\
=& C_1(\E \lambda_{t}^{2} - (\E \lambda_{t})^{2}) + C_1(\E \lambda_{t}- \lambda_{t}^{*})^{2} + C_2 (\E \lambda_{t} - \lambda_{t}^{*}).
\end{aligned}
\end{equation}
where
$$C_1 = \frac{N}{2} \left( \sum_{i=1}^N \frac{1}{\beta_i} + ( \sum_{i=1}^N  \frac{1}{\beta_i} )^{2} \right),$$
and
$$ C_2 = \sum_{i=1}^N \frac{\alpha_i}{\beta_i} \left( \sum_{i=1}^N \frac{1}{\beta_i} -1 \right) $$ are some constant coefficients that do not evolve with time $t$.

From \eqref{regretcase2}, we observe that the regret consists of variance, bias and squared bias of $\lambda_t$, respectively. This suggests that it is preferable to have a pricing strategy that achieves both small variance and bias, or can tradeoff  two.





\subsection{Online learning procedure}

\begin{figure}
	\centering
\begin{tikzpicture}[node distance = 4cm, auto]

\node [block] (init) {index $t = 1$, history observation $\mathcal{H}_{t-1}$ before time $t$};
\node [block, below of=init] (identify) {At time $t$, utility updates parameter estimation $\hat{\gamma}_1^t$  and $\hat{\gamma}_1^t$  from Algo. \ref{lin1}, and update the price signal as:
$$
\lambda_t = \frac{Y^*d_t + \hat{\gamma}_2^t}{\hat{\gamma}_1^t + N},
$$ and
$\lambda_t$ is broadcasted to users.};
\node [block, below of=identify, node distance = 3.2cm] (response) {User $i$ responds to the reward signal and reveals its consumption as:
	$$
	\hat{x}_{i}^{t} = \frac{N \lambda_t - \alpha_i}{\beta_i} + \epsilon_i^t,
	$$};
\node [smallblock, below of=response, node distance = 2.2cm] (evaluate) {$t = t + 1$};
\node [smallblock, left of=identify, node distance=4cm] (update) {update history $\mathcal{H}_{t-1}$};
\node [decision, below of=evaluate] (decide) {$t = T$?};
\node [smallblock, below of=decide, node distance=2.2cm] (stop) {stop};
\path [line] (init) -- (identify);
\path [line] (identify) -- (response);
\path [line] (response) -- (evaluate);
\path [line] (evaluate) -- (decide);
\path [line] (decide) -| node [near start] {No} (update);
\path [line] (update) |- (identify);
\path [line] (decide) -- node {Yes}(stop);


\end{tikzpicture}
\caption{Online learning algorithm for pricing.} \label{fig_algo1}
\end{figure}
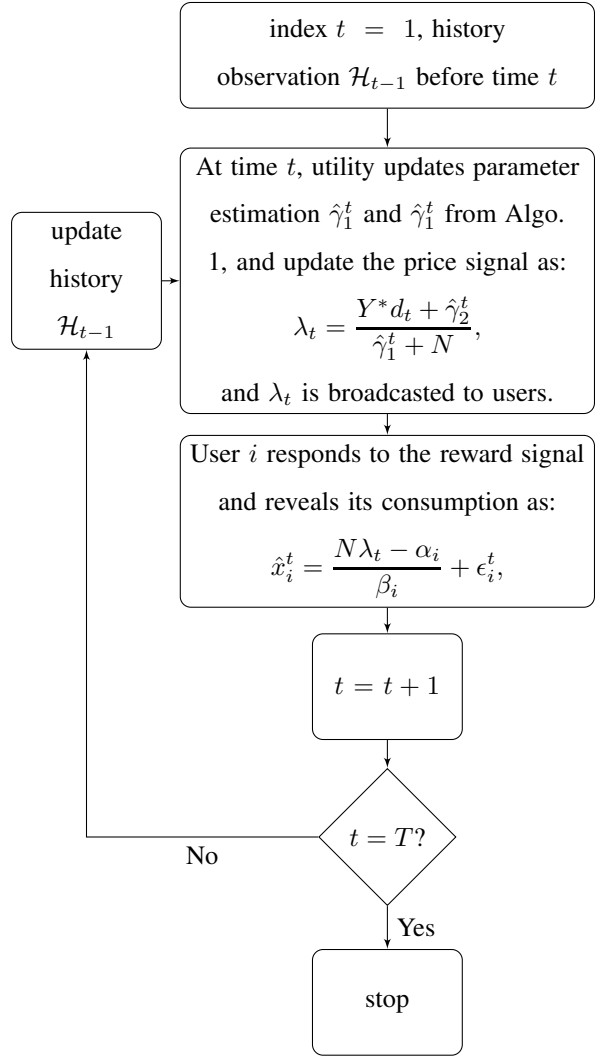

The online learning strategy for optimal prices is presented in Fig. \ref{fig_algo1}. We denote $ \gamma_1 \triangleq \sum_i \frac{1}{\beta_i}$ and $\gamma_2 \triangleq -\sum_i \frac{\alpha_i}{\beta_i}$.
Since at each time period the users' response is linear in the price signal, we propose to estimate the unknown parameters through linear regression based on history using least squares. The parameters to estimate are: $\gamma_1, \gamma_2 $. 
The estimation of $ \gamma_1$ and $ \gamma_2$ at time $t$ are denoted by $\hat{\gamma}_1^t$ and $\hat{\gamma}_2^t$, respectively.

The online learning algorithm shown in Fig. \ref{fig_algo1} consists of several steps: first, the utility collects history observation until the current time point; then it adopts the decision from Algo. \ref{lin1} (linear regression) to design the price signal and broadcast it to the end users. Finally, the users' consumption based on this price signal is reported back to the utility and the process repeats.

The core of the online algorithm is in Algo. \ref{lin1}, which determines the performance of the algorithm. In Algo. \ref{lin1} we estimate the unknown parameters $\gamma_1$ and $\gamma_2$ by the least square method\footnote{When $t$ = 1 and 2, since there are no enough data points to derive least square estimator from the linear regression model, we impose a prior information on the parameters and estimation is done as ridge regression.}. What is more, the whole procedure is done iteratively by adding more samples into the model, which means that as $t$ gets bigger, we accumulate more samples to train the model and iterate. We thus call Algo.\ref{lin1} \emph{iterative linear regression}. Its details are illustrated in the next subsection. 

The performance of the algorithm in Fig. \ref{fig_algo1} is evaluated by regret. Recall that Theorem \ref{maintheorem} states there exists an algorithm that the regret $R$ is growing logarithmically with time horizon $T$. We use Theorem \ref{maintheorem2} to show that the proposed algorithm in Fig. \ref{fig_algo1} achieves this rate.

Moreover, we adopt the following asymptotic bound notations. These notations facilitate the analysis to compare the orders of quantities of interest~\cite{algo}.
\begin{definition}\label{def_Theta}
	$f(n) = \Theta (g(n)) $ means there are positive constants $k_1$, $k_2$, and $n_0$, such that $0 \leq k_1g(n) \leq f(n) \leq k_2g(n), \forall n \geq n_0$.
\end{definition}
\begin{definition}\label{def_bigO}
	$f(n) = O (g(n)) $ means there are positive constants $k$, and $n_0$, such that $ f(n) \leq k_2g(n), \forall n \geq n_0$.
\end{definition}

With these notations, we state Theorem \ref{maintheorem2} as below.

\begin{theorem}\label{maintheorem2}
	The algorithm presented in Fig. \ref{fig_algo1} achieves $\Theta(\log T)$-regret.
\end{theorem}

In addition, recall that, if $f(n) = \Theta(g(n))$ then $f(n) = O(g(n))$, the proposed algorithm at the same time achieves a regret of $O(\log T)$ as well. Therefore, proving Theorem \ref{maintheorem2} infers Theorem \ref{maintheorem}. The detailed proof of Theorem 2 can be found at Appendix. Here we provide a sketch of the proof to Theorem 2.
\begin{proof}[Sketch of proof to Theorem \ref{maintheorem2} ]
	We want to show that at each time $t$, the gap between the offline optimal strategy and online strategy is bounded by $\frac{1}{t}$. Then following the fact that $\sum_t \frac{1}{t} = \log T$ for $T$ time slots, the regret is $\log T$. To show that the gap is $\frac{1}{t}$, from \eqref{regretcase2} it suffices to show that the variance and the bias (with a multiplicative term $C_1$ and $C_2$) of the price estimate $\lambda_t$ is decaying at a rate of $\frac{1}{t}$ at each time slot $t$. In the detailed proof shown in the Appendix, we show that this rate indeed holds.
\end{proof}

\subsection{Iterative linear regression}

As can be seen from Fig. \ref{fig_algo1}, 
the performance of the algorithm in Fig. \ref{fig_algo1} largely depends on the estimator of $\gamma_1$ and $\gamma_2$ obtained by the least square from Algo. \ref{lin1}, since they influence variance and bias of $\lambda_t$ which determine the regret. Algo.\ref{lin1} is shown with one step iteration.


	\begin{algorithm}[h]
		\caption{Iterative linear regression (one step)}
		\label{lin1}
		\begin{algorithmic}[1]
			\State \textbf{Input}: History $\mathcal{H}_{t-1}$ which includes price history ${\lambda}_s$ and response history sequence $\sum_i{\hat{x}_i^s}$ for  $s = 1, 2, ..., t-1$ and $t$ is the current time stamp.
			\State The linear regression model is:
			\begin{equation}\label{linreg2}
			\sum_i \hat{x}_i^s = N {\lambda}_s {\gamma}_1 + {\gamma}_2 + \sum_i \epsilon_i^s , s \in \{1,2,...,t-1\}.
			\end{equation}
			\State Do a least square estimate in the linear regression model of  $\sum_i{\hat{x}_i^t}$ on $N{\lambda}_t$ plus an intercept. 
			\State \textbf{Output}: Least square estimate $\hat{\gamma}_1^t$ and $\hat{\gamma}_2^t$.
		\end{algorithmic}
	\end{algorithm}

There are a few points to note. First, the estimators from Algo. \ref{lin1} are consistent but not unbiased. This is because that the online prices $\lambda_t$ that we generate based on the estimates from Algo. \ref{lin1} is correlated with past noise. These $\lambda_t$'s are again fed into Algo. \ref{lin1} as regressors to train the linear regression model. However, the bias is of lower orders compared to the order of the estimator and is decaying with the number of time periods, we approximate the estimators as well as their variance and expectation by assuming that they are unbiased estimators. This approximation is validated by simulation results.

Second, one thing that might hinder one from using iterated linear regression in online learning algorithms is that it may happen that the observations are not exploring the linear regression model. For example $\lambda_t$ is the same for all $t \in \{1,2,...,s\}$. If so, the linear regression is not efficient because the two regressors ($\lambda_t$ and the intercept) is the same which renders infinite variance in the estimator. This problem is addressed in \cite{Keskin2014, KhezeliEtAl2016}.

However, with the assumption that $d_t$ is randomly spread between $[d_{min}, d_{max} ]$, \eqref{optimallambda} suggests that $\lambda_t$ will also be much different across different time $t$. Then we can sufficiently explore the structure of the linear regression model and that the estimation from least square is efficient. In the simulation, we show that even some of the $d_t$'s are similar, the exploration of the linear regression model is still effective, such that the regret is still in order $\log T$. 

	\section{Simulation Results}\label{simulation}
\subsection{Parameter set up}

In the simulation, we generate a random price at $t_1 = 1$ to start the online algorithm. For practical purposes, we introduce ridge regression estimators in the linear regression model \cite{ridge}. The influence of $\lambda$ on the estimator will decay fast as more samples are included into the model. For interested readers, the computation of ridge regression estimator can be found in the Appendix.

The parameters of the system are shown in Table \ref{table2}, where $d_t$ are the normalized demand reduction and $c$ is a constant bounded away from zero. The iteration of Monte Carlo simulation for the responses of users is set to be 1000 and the regularization parameter in the ridge regression is set to be 0.001. In addition, there are totally 100 different users.
\begin{table}[!ht]
	\renewcommand{\arraystretch}{1.3}
	\caption{Parameters. Intervals indicate the uniform distribution. $c$ is some constant. We simulate two sets of parameters and compare the results.}
	\centering
	\begin{tabular}{|c|c|c|c|}
		\hline
		\bfseries $\alpha$ & \bfseries $\alpha_i$ & \bfseries  $\beta_i$ & \bfseries  $d_t$ \\
		\hline
		$c \max \{d_t\}$  & [1,2] &  [4,8] & [3,6]\\
		\hline
		$c \max \{d_t\}$  & [1,3] &  [3,10] & [2,5]\\
		\hline
	\end{tabular}
	\label{table2}
\end{table}

\subsection{Pricing strategy}
\begin{figure}[!ht]
		\centering
		\begin{subfigure}[b]{\textwidth}
			\centering
			\includegraphics[width=0.8\linewidth]{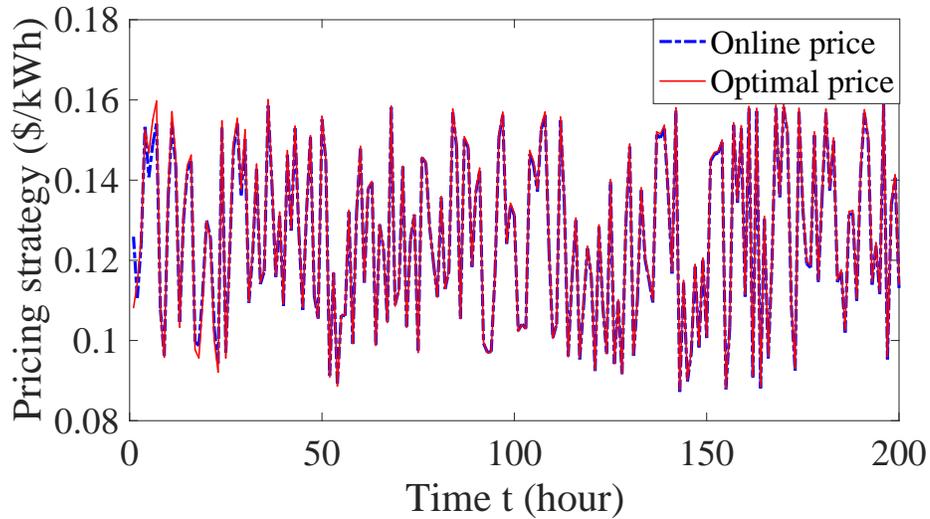}
			\caption{Optimal pricing and online pricing.}
			\label{price_comp}
		\end{subfigure}
		\begin{subfigure}[b]{\textwidth}
			\centering
			\includegraphics[width=0.8\linewidth]{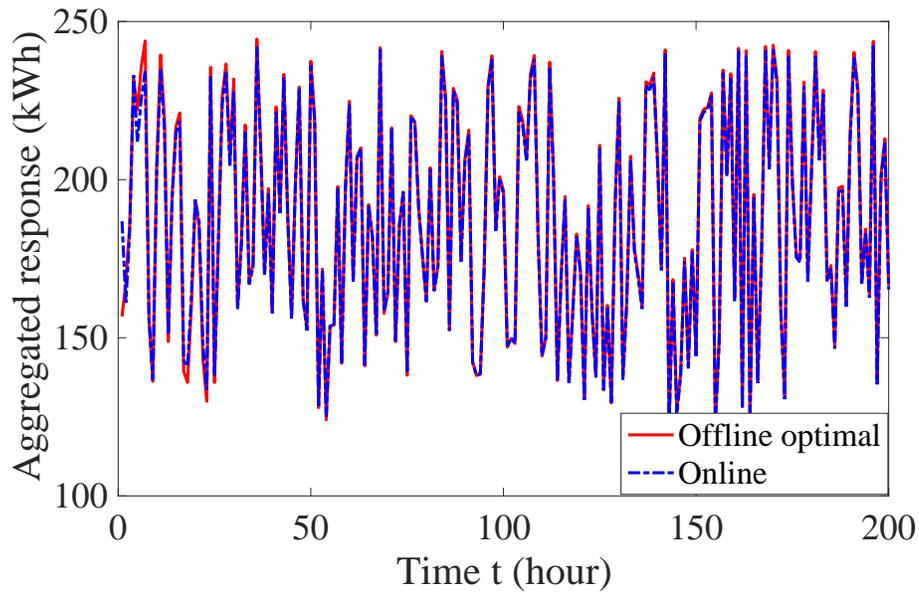}
			\caption{Optimal response and online response (aggregated).}
			\label{aggregated_response}
		\end{subfigure}

		\caption{Comparison between optimal pricing (aggregated response) and online pricing.}
	\end{figure}

We first validate the online pricing strategy discussed in Section \ref{online}. The comparisons between the optimal pricing and online pricing are shown in Fig.\ref{price_comp} and Fig. \ref{aggregated_response}. As can be seen from Fig.\ref{price_comp} and Fig. \ref{aggregated_response}, at the first few time steps, the online algorithm is still learning the parameters given very few sample points, and the online price deviates a lot from the optimal pricing. This deviation in price also drives the aggregated response to be suboptimal. However, the estimation accuracy is greatly improved given just a few more sample points, and as we can see from Fig.\ref{price_comp} and Fig. \ref{aggregated_response}, the online price tracks the optimal price after less than 50 time steps, so does the aggregated response. This means that after a few time steps, the aggregated response is able to achieve the demand response requirement by the utility, even the utility does not know the user's parameters in advance.

\begin{figure}[!ht]
	\centering
	\includegraphics[width = 0.8\columnwidth]{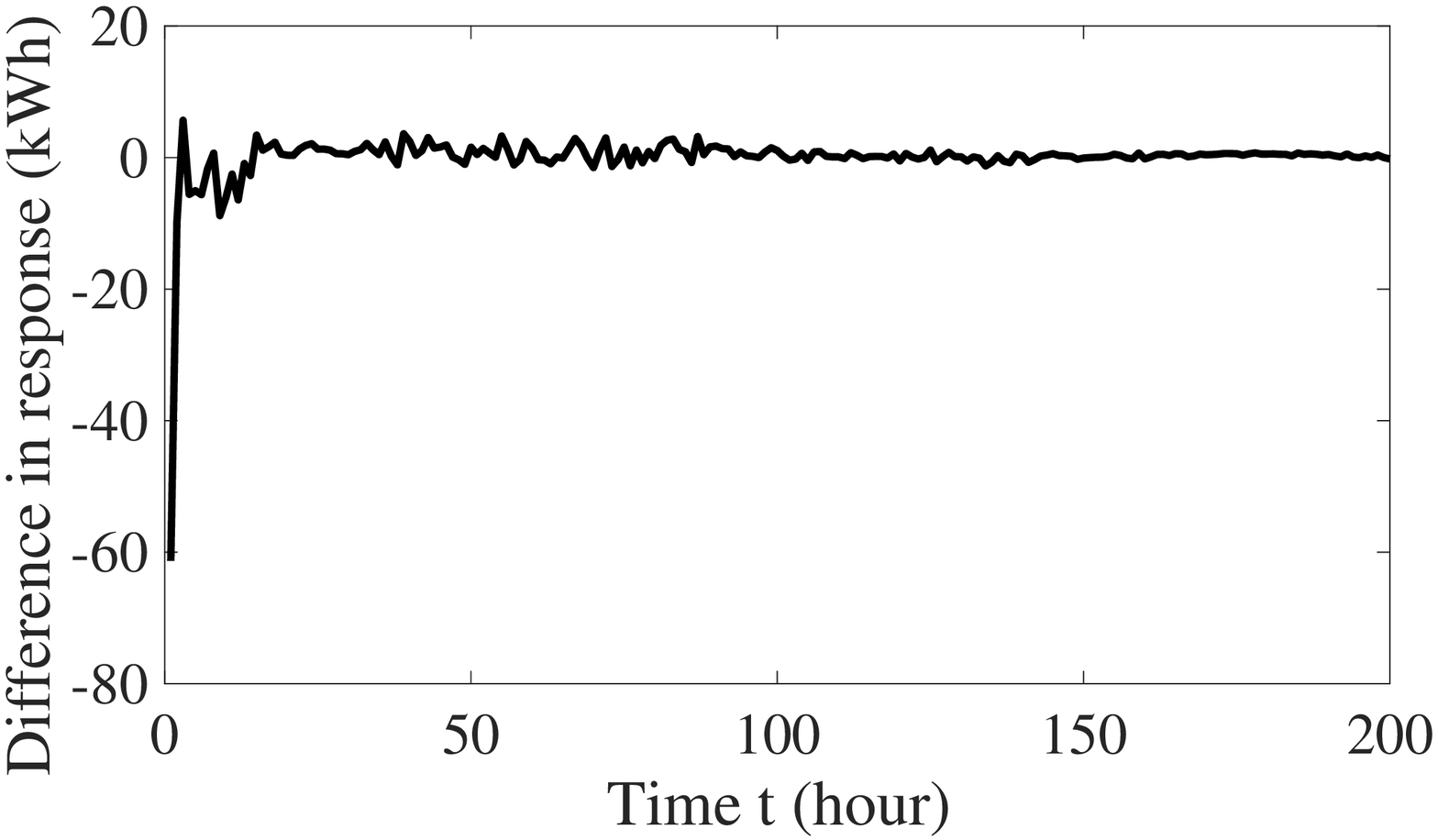}
	\caption{Difference between the aggregated response induced by online pricing and optimal pricing.}
	\label{diff}
	\end{figure}

The comparison is clearer when we explicitly compute the difference between aggregated online response and optimal response, shown in Fig. \ref{diff}. As can be seen from Fig. \ref{diff}, the aggregated response $\sum_i \hat{x}_i^t(\lambda_t)$ induced by the proposed online pricing strategy $\lambda_t$ is approaching the optimal response $\sum_i \hat{x}_i^t(\lambda_t^*)$. The difference between the two responses is diminishing and becomes zero as time goes on.

		\begin{figure}[!ht]
			\centering
			\begin{subfigure}[b]{\textwidth}
				\centering
				\includegraphics[width=0.6\linewidth]{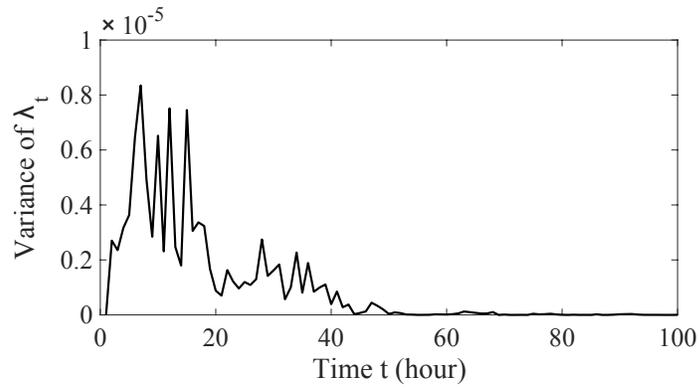}
				\caption{Variance of $\lambda_t$}
				\label{var_lambda}
			\end{subfigure}
			\begin{subfigure}[b]{\textwidth}
				\centering
				\includegraphics[width=0.6\linewidth]{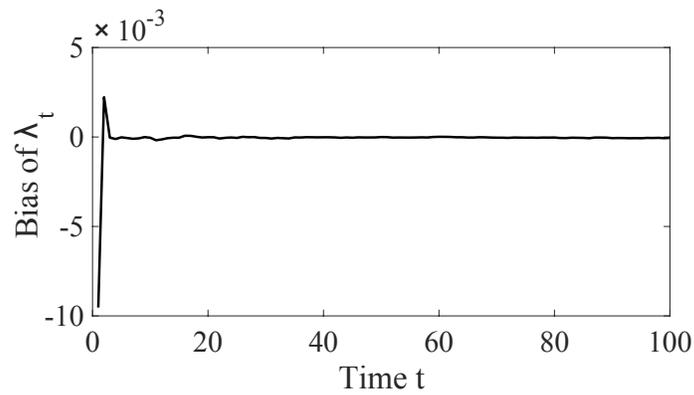}
				\caption{Bias of $\lambda_t$}
				\label{bias_lambda}
			\end{subfigure}

			\caption{Variance and bias of online price $\lambda_t$.}
		\end{figure}

	The variance and the bias of online price $\lambda_t$ is shown in Fig. \ref{var_lambda} and Fig. \ref{bias_lambda}. As we can see, at first the bias is huge because $\lambda_1$ is randomly initialized, so it is not accurate. As long as the learning procedure begins, the bias drops drastically. The variance of $\lambda_t$ has a much slower decay rate of $\frac{1}{t}$, as has been discussed in Section \ref{online}. Comparing the orders, the squared bias is much smaller than the variance. Therefore the tradeoff between variance and squared bias is dominated by variance.

\subsection{Regret Analysis}
\begin{figure}[ht!]
	\centering
	\begin{subfigure}[t]{\textwidth}
		\centering
		\includegraphics[height=2.4in]{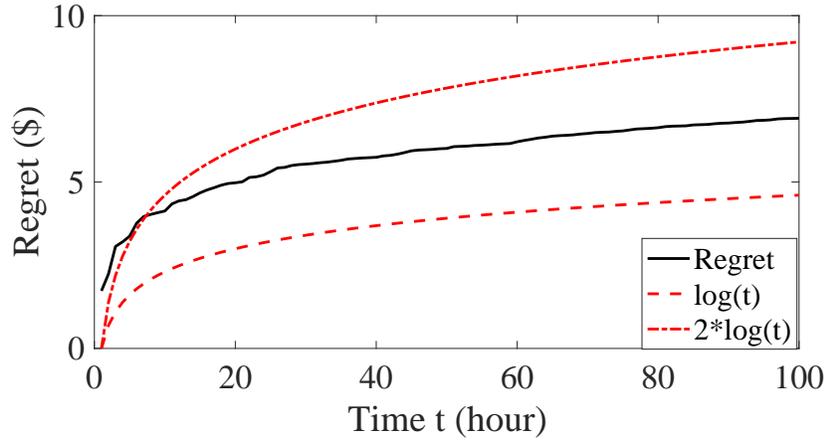}
		\caption{$\alpha_i \sim [1,2]$, $\beta_i \sim [4,8]$, $d_t \sim [3,6]$.}
	\end{subfigure}%
	\quad
	\begin{subfigure}[t]{\textwidth}
		\centering
		\includegraphics[height=2.4in]{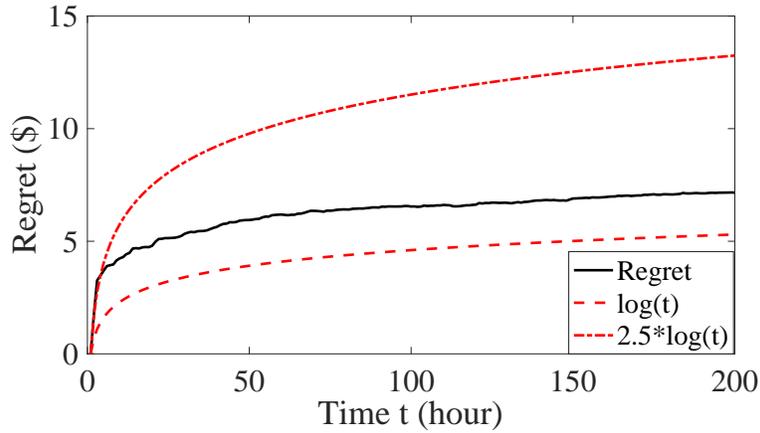}
		\caption{$\alpha_i \sim [1,3]$, $\beta_i \sim [3,10]$, $d_t \sim [2,5]$.}
	\end{subfigure}
	\caption{Regret over time with different parameters.}
	\label{cumregret2}
\end{figure}

We then analyze the regret for the online pricing strategy. The gap $R_t$ between online cost and offline cost is shown in Fig. \ref{eg1} in Section \ref{mainresults} and the regret is shown in Fig. \ref{cumregret2}.

From Fig. \ref{eg1}, we observe that the gap decays fast with time with a rate of $\frac{1}{t}$, which means the cost obtained from online price signals $\lambda_t$ approaches the true cost to the system quickly with time. The regret is the sum of the gaps during all time periods and is in order $\log T$ as shown in Fig. \ref{cumregret2}. As can be seen from Fig. \ref{cumregret2}, the regret is within some bounds of $\log t$, which validates Theorem \ref{maintheorem2}. What is more, from Fig.\ref{cumregret2}, we find that the algorithm works for different sets of parameters.

\subsection{Performance of the algorithm subject to same consecutive requirement over time}

As has been pointed out in \cite{KhezeliEtAl2016}, an iterated linear regression may not well explore the model. This results in a larger variance in the estimator $\hat{\gamma}_1$ and $\hat{\gamma}_2$, and thus deteriorates the estimation of price, i.e., $\lambda_t$, which may lead to a regret worse than $O(\log T)$.

With a more careful scrutiny, we find that this situation only happens when the regressors $\lambda_t$'s are highly correlated. We have argued in Section \ref{online} that this situation is avoided with the assumption that $d_t$'s are different across time, or there are sufficiently many $d_t$'s that are not the same to each other.

\begin{figure}[!ht]

	\centering
	\includegraphics[width = 0.8\columnwidth]{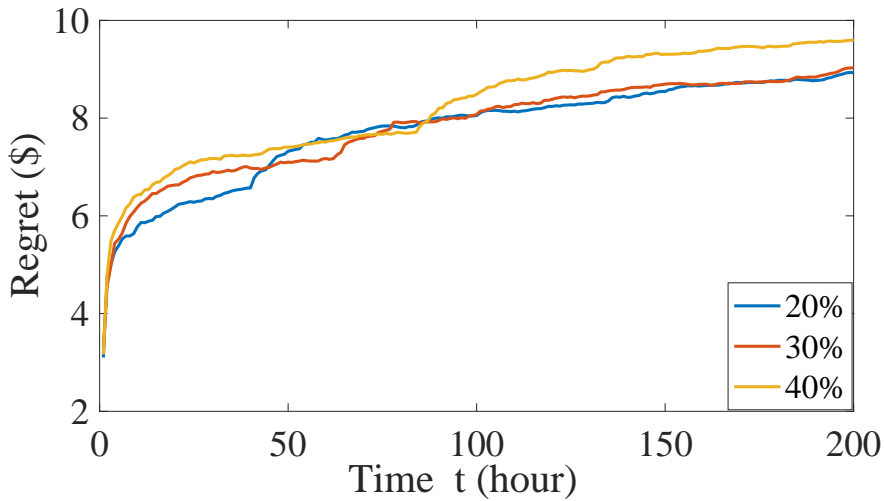}
	\caption{Regret $R$ over time where 20 \%, 30\% and 40\% of $d_t$'s are the same. Cost function is quadratic and linear.}
	\label{regret_diff_same_dt}
\end{figure}

To state this clearer, we first set up a sequence of $d_t$'s in which 20 \%, 30\% and 40\% of them are identical. The result is shown in Fig. \ref{regret_diff_same_dt}.



From Fig. \ref{regret_diff_same_dt}, we see that the regret is still sub-linear, i.e., $\Theta(\log t)$, in time $t$. This suggests that even with some portion of same $d_t$'s in the system, the proposed algorithm still works effectively.

\begin{figure}[!ht]

	\centering
	\includegraphics[width = 0.8\columnwidth]{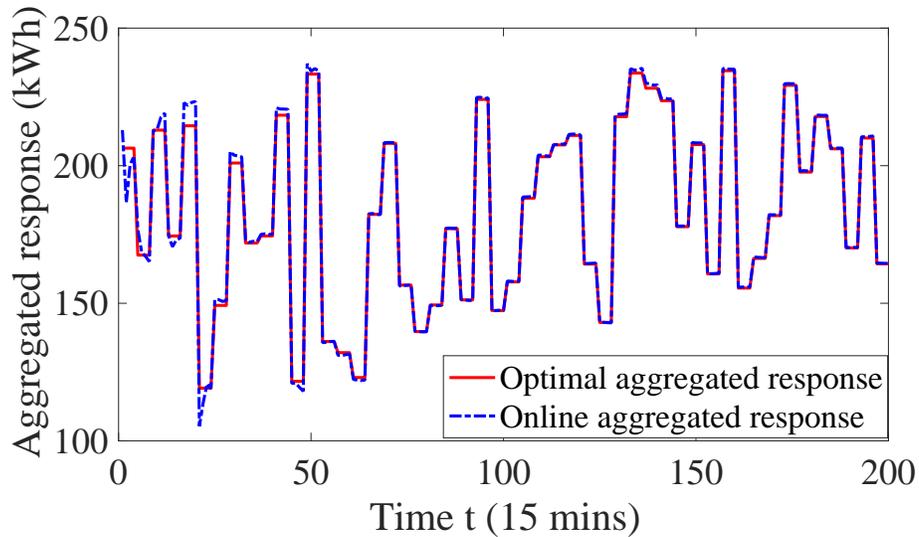}
	\caption{Comparison of the optimal aggregated response and online aggregated response, where $d_t$ is the same every four time slots.}
	\label{review_3_13}
\end{figure}

\begin{figure}[!ht]

	\centering
	\includegraphics[width = 0.8\columnwidth]{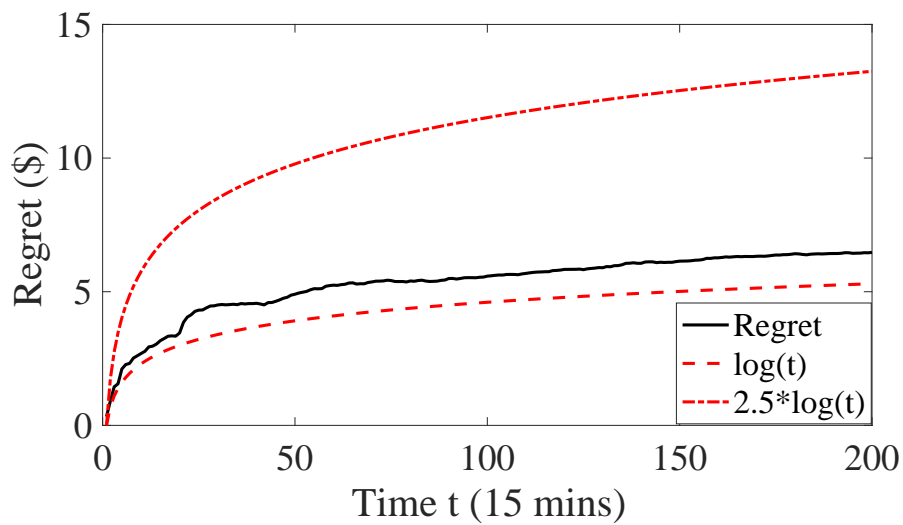}
	\caption{Regret $R$ where $d_t$ is the same every four time slots.}
	\label{review_3_13_regret}
\end{figure}

We also test against a case where each time period is short, i.e., $t$ represent fifteen minute-level observation. Suppose that the demand response requirement only get changed every hour, which means that $d_t$ remains the same for every four time slots. In total, 25\% of all the $d_t$'s are the same. The results are shown in Fig. \ref{review_3_13} and Fig. \ref{review_3_13_regret}.

As can be seen from Fig.\ref{review_3_13}, the response now is smoother than that shown in Fig. \ref{aggregated_response}, since the observation is more frequent. The online aggregated response is again able to achieve the optimal response, after roughly 50 time slots.

The regret is shown in Fig. \ref{review_3_13_regret}. It is clearly within the bounds of $\log T$, which validates the statement that the algorithm remains effective with some amount of repetitive pattern in $d_t$.

	\section{Conclusion}\label{conclusion}

We studied an online demand response strategy for a utility company. We developed an optimization framework for the demand response utility to minimize the system cost, including users' costs and utility's cost. We considered a realistic setting where the utility does not know the cost functions of users. The utility thus makes sequential pricing decisions to learn the cost functions of users as well as minimizing the cost at the same time. We designed an online learning algorithm referred to as iterated linear regression to assist this pricing strategy. Compared with the offline optimal solution, our online algorithm achieves a sublinear regret in the length of time period, i.e., $\Theta (\log T)$. We validated the theoretic results through numerical simulations. Further simulation results showed that the algorithm is efficient even when the consumption requirement is the same across multiple time periods. The proposed framework and the online algorithm sheds light on more complicated settings where the utility selectively calls users into the demand response, and is left as extended work to this paper.

\begin{appendix}

\subsection{Preliminaries for proving Theorems \ref{maintheorem} and \ref{maintheorem2}}
Note that Theorem 2 directly implies Theorem 1 since it directly characterizes the performance of the online learning algorithm. We separate the proof of Theorem \ref{maintheorem2} into two steps, where we first show it for the case $\alpha_i=0$ for every $i$,  then for the case of general $\alpha_i$'s.


Before proving Theorem \ref{maintheorem2}, we need to introduce a few definitions that help the proof. We use asymptotic bounds to analyze the interactions between quantities in the following analysis, since many random quantities are involved and an exact inference may not be possible. The asymptotic bound notation that we use in this paper is $\Theta(\cdot)$ and $O(\cdot)$ as defined in Section \ref{offline}. These bounds quantify the orders of the terms involved in computation and ignore constant factors, which simplify computation.

In addition, we derive some useful operations for these notations.
\begin{remark}\label{remark2}
If for some positive functions $f_1(n)$, $g_1(n)$, $f_2(n)$,$g_2(n)$, $f_1(n) = \Theta (g_1(n))$ and $f_2(n) = \Theta (g_2(n))$ with same $n_0$, then $\frac{f_1(n)}{f_2(n)} = \Theta(\frac{g_1(n)}{g_2(n)})$. In addition, $f_1(n)^2 = \Theta (g_(n)^2)$ and $f_1(n)f_2(n) = \Theta(g_1(n)g_2(n))$.
\end{remark}
%
%
%
%

Using the asymptotic bound notation we greatly simplify the comparison between different quantities in each of the equations presented earlier. For example, we can write $\hat{\gamma}_1^t$ as a result from least square estimation in the linear regression model presented in \eqref{linreg2}, where $\alpha_i = 0$:
$
\hat{\gamma}_1^t = \frac{\sum_{s=0}^{t-1}  (\sum_i x_i^s) N \lambda_s}{\sum_{s=0}^{t-1} N^2 \lambda_s^2}.
$
Its variance is expressed as:
$
\Var (\hat{\gamma}_1^t) = \frac{N}{\sum_{s=0}^{t-1} N^2 \lambda_s^2}.
$

If we know the order of $\lambda_t$ in terms of $\Theta (\cdot)$, etc., then the order of the estimator and its variance can be obtained, which helps to determine the order of the regret $R$.

\begin{remark}\label{remark1}
	We assume that $d_t$ is bounded away from zero. What is more, assume that $\alpha_i$ and $\beta_i$ are also bounded away from zero and take some finite values. Along with \eqref{optimal2}, \eqref{Ystar} and assume that the revenue $\alpha$ is bounded below from zero, we obtain that $Y^* = \Theta (N)$, which leads to the fact that $\lambda_t = \Theta( \frac{1}{N})$ with high probability.
\end{remark}

%

\subsection{Proof of Theorem \ref{maintheorem2} with $\alpha_i=0$}
\begin{proof}
	With Remark \ref{remark1}, we obtain the least square estimator for $\gamma_1$ at time $t$ as the following:
	\begin{equation}\label{varcase1}
	\begin{aligned}
	\Var (\hat{\gamma}_1^t) & = \frac{N}{\sum_{s=0}^{t-1} N^2 \lambda_s^2} \\
	& = \frac{1}{N} \frac{1}{\sum_{s=0}^{t-1} \lambda_s^2} \\
	& \overset{(a)} = \frac{1}{N} \frac{1}{\lambda_0^2 + (t-1) \left[ \Theta \left( \frac{(Y^*)}{N^2} \right)  \right]^2} \\
	& \overset{(b)} = \frac{1}{N} \frac{1}{\lambda_0^2 + (t-1) \Theta \left( \frac{(Y^*)^2}{N^4} \right)}\\
	& \overset{(c)} = \Theta \left( \frac{N^3}{t (Y^*)^2} \right)\\
	& \overset{(d)} = \Theta \left( \frac{N}{t} \right),
	\end{aligned}
	\end{equation}
	where (a) follows from the fact that since $\sum_i \frac{1}{\beta_i} = \Theta (N)$ then from \eqref{optimallambda} we know that $\lambda_s = \Theta(\frac{Y}{N^2})$. Equality (b) follows from Remark \ref{remark2} and (c) is the direct result from further simplification using $\Theta$ notation. Finally (d) comes from the fact that $Y^* = \Theta(N)$.

	Using the same simplification from asymptotic bounds and based on the result in \eqref{varcase1} , we can rewrite the gap $R_t$ in \eqref{regrett} as
	\begin{equation}
	\begin{aligned}
	R_t  =& C_1 \left[ \E \lambda_{t}^{2} - (\E \lambda_{t})^{2} + (\E \lambda_{t}- \lambda_{t}^{*})^{2}  \right] \\
	\overset{(a)} =& C_1
	 \left[ \Var \left( \frac{Y^* d_t}{N \hat{\gamma}_1^t + N} \right) \right. \\
	 & \left. + \left( \E \left( \frac{Y^* d_t}{N \hat{\gamma}_1^t + N} \right) - \frac{Y^* d_t}{N ( \sum_{i \in \mathcal{N} }  \frac{1}{\beta_i} ) + N} \right)^2 \right] \\
	\overset{(b)} =& C_1
	 \left[ \frac{(Y^*)^2 d_t^2 N^2 \Var \hat{\gamma}_1^t}{ (N \sum_{i \in \mathcal{N} }  \frac{1}{\beta_i} + N)^4 } + \frac{(Y^*)^2 d_t^2 N^4 (\Var \hat{\gamma}_1^t)^2}{(N \sum_{i \in \mathcal{N} }  \frac{1}{\beta_i}t+ N)^3}  \right] \\
	\overset{(c)} =& C_1
	 \left[ \frac{(Y^*)^2 d_t^2 N^2 \Var \hat{\gamma}_1^t}{ \Theta (N^8) } + \frac{(Y^*)^2 d_t^2 N^4 (\Var \hat{\gamma}_1^t)^2}{\Theta (N^{12})}  \right]  \\
	\overset{(d)} =& C_1
	 \left[ \Theta(\frac{1}{t N^3}) + \Theta(\frac{1}{Y^2 N^2 t^2}) \right],
	\end{aligned}
	\end{equation}
	where (a) follows according to the result in \eqref{optimallambda}, (b) follows from the second order approximation of the expectation of an inverse random variable \cite{ratio}. Then using Remark \ref{remark1} and Remark \ref{remark2}, we arrive at equality (c) and (d).

	Given that $Y^*= \Theta (N)$ and $C_1 = (\frac{N}{2} \sum_{i =1}^N  \frac{1}{\beta_i})  + \frac{N}{2} ( \sum_{i =1}^N  \frac{1}{\beta_i} )^{2} = \Theta(N^3)$, we have
	\begin{equation}\label{regret1_final}
	\begin{aligned}
	R_t =& \Theta (N^3) \left[\Theta( \frac{1}{t N^3}) + \Theta (\frac{1}{Y^2 N^2 t^2}) \right ] \\
	=& \Theta (\frac{1}{t} + \frac{1}{N t^2}).
	\end{aligned}
	\end{equation}

	Since $\frac{1}{Nt^2} = O(\frac{1}{t})$, $ R_t = \Theta(\frac{1}{t})$. This indicates that the gap $R_t$ forms a harmonic series with respect to time, which means that the cumulative regret $R = \sum_t R_t$ is
	$\Theta (\log T)$. This also implies that $R = O(\log T)$.
\end{proof}

\subsection{Proof of Theorem \ref{maintheorem2} with $\alpha_i \neq 0$}
\begin{proof}
	First write the linear regression model in \eqref{linreg2} in a more compact form as:
	\begin{equation}
		\begin{aligned}
		Z_t = \begin{bmatrix}
		\gamma_1 & \gamma_0
		\end{bmatrix}
		\begin{bmatrix}
		N {\lambda}_t \\
		1
		\end{bmatrix} + \sum_{i \in \mathcal{N}} \epsilon_i^t,
		\end{aligned}
	\end{equation}
	where $Z _t\triangleq \sum_i \hat{x}_i^t$,  $\gamma_1 = \sum_i \frac{1}{\beta_i}$, $\gamma_0 = - \sum_i \frac{\alpha_i}{\beta_i}$. The white noise in the model is $\sum_i \epsilon_i^t \sim \mathcal{N}(0, \sigma^2)$ and $\sigma^2 = N$.

	The estimator at time $t$ is obtained from least squares in the following form:
	\begin{equation}
	\begin{aligned}
	\begin{bmatrix}
	\hat{\gamma}_1 \\
	\hat{\gamma}_0
	\end{bmatrix}
	= (\mathbf{X}^\top \mathbf{X})^{-1} \mathbf{X}^\top \mathbf{Z},
	\end{aligned}
	\end{equation}
	where $\mathbf{X} = \begin{bmatrix}
	N \bm{\lambda}_t \quad \mathbf{1}
	\end{bmatrix}$, $\bm{\lambda}_t =  [\lambda_1, \lambda_2,  \dots, \lambda_t]^\top$ is a column vector and $\mathbf{1}$ is an all one column vector, $\mathbf{Z} = [Z_1, Z_2, ..., Z_t]^\top$.

	We calculate
	\begin{equation}\label{var2}
	\begin{aligned}
	\Var
	\begin{bmatrix}
	\hat{\gamma}_1^t \\
	\hat{\gamma}_0^t
	\end{bmatrix}
	= &
	\begin{bmatrix}
	\Var \hat{\gamma}_1^t  & \Cov (\hat{\gamma}_1^t , \hat{\gamma}_0^t ) \\
	\Cov (\hat{\gamma}_0^t , \hat{\gamma}_1^t ) & \Var \hat{\gamma}_0^t
	\end{bmatrix}
	\\
	=& \left( \mathbf{X}^\top \mathbf{X} \right)^{-1} \sigma^2 \\
	=& \left(
	\begin{bmatrix}
	N \bm{\lambda}_t^\top \\
	\mathbf{1}^\top
	\end{bmatrix}
	\begin{bmatrix}
	N \bm{\lambda}_t &\mathbf{1}
	\end{bmatrix}
	\right)^{-1} \sigma^2 \\
	=&
	\begin{bmatrix}
	N \bm{\lambda}_t^\top N \bm{\lambda}_t & N \bm{\lambda}_t^\top \mathbf{1} \\
	\mathbf{1}^\top N \bm{\lambda}_t &\mathbf{1}^\top \mathbf{1}
	\end{bmatrix} ^{-1} \sigma^2 \\
	\overset{(a)}=&
	\frac{N}{t\sum_s(N\lambda_s)^2 - (\sum_s N\lambda_s)^2} \\
	& \cdot \begin{bmatrix}
	t & - \sum_s N\lambda_s\\
	- \sum_s N\lambda_s & \sum_s(N\lambda_s)^2\\
	\end{bmatrix} ,\\
	\end{aligned}
	\end{equation}
	where equality (a) follows from the matrix inversion and that $\sigma^2 = N$.

	Since we argue that $\lambda_t$'s are not the same across different $t$ (given that $d_t$'s are not the same across time), we obtain that $t\sum_s(N\lambda_s)^2 - (\sum_s N\lambda_s)^2$ is bounded below from 0, which validates the least squares estimator from linear regression. What is more, since each $\lambda_t = \Theta(\frac{1}{N})$ according to \eqref{optimal2} and \eqref{optimallambda}, we know that $t\sum_s(N\lambda_s)^2 - (\sum_s N\lambda_s)^2 = \Theta (\frac{1}{t^2})$. Plug this result back in \eqref{var2}, we obtain that the estimators at any time $t$ have the following asymptotic bounds on their variances and covariance:
	\begin{equation}\label{var2final}
		\Var \hat{\gamma}_1^t , \Var \hat{\gamma}_0^t ,  \Cov (\hat{\gamma}_1^t , \hat{\gamma}_0^t ) = \Theta(\frac{N}{t}).
	\end{equation}

	Therefore, the bias of $\lambda_t$ is calculated as:
	\begin{equation}\label{bias2}
	\begin{aligned}
	\text{bias} \lambda_t
	=& \E \lambda_t - \lambda_t \\
	=& \E \left( \frac{Y^* d_t + \hat{\gamma}_{0}^t }{N \hat{\gamma}_{1}^t  + N} \right) -  \frac{Y^* d_t + \gamma_{0}}{N \gamma_{1} + N} \\
	\overset{(a)}=& - \frac{\mathrm{Cov} (\hat{\gamma}_{0}, N \hat{\gamma}_{1}^t )}{\left( \E (N \hat{\gamma}_{1}^t  + N) \right)^2}
	+ \frac{\Var(N \hat{\gamma}_{1}^t ) \E(Y^* d_t + \hat{\gamma}_{0}^t )}{(N \gamma_{1} + N)^3} \\
	\overset{(b)}=&  \Theta(\frac{N \cdot \frac{N}{t}}{N^4}) +  \Theta (\frac{N^2 \cdot \frac{N}{t} \cdot N}{N^6}) \\
	=&  \Theta (\frac{1}{N^2 t}),
	\end{aligned}
	\end{equation}
	where equality (a) comes from second order approximation of the expectation of a ratio between two random variables as in \cite{ratio} and equality (b) follows from Remark \ref{remark1} and Remark \ref{remark2}.

	Using the same analysis as when $\alpha_i = 0$, and based on \eqref{var2final} and \eqref{bias2}, we obtain the gap $R_t$ as:
	\begin{equation}
	\begin{aligned}
	R_t =& \frac{N}{2} \left( \sum_{i =1}^N  \frac{1}{\beta_i} + (  \sum_{i =1}^N   \frac{1}{\beta_i} )^{2} \right)
	\left[ \E \lambda_{t}^{2} - (\E \lambda_{t})^{2} + (\E \lambda_{t}- \lambda_{t}^{*})^{2}  \right] \\
	&+  \sum_{i =1}^N  \frac{\alpha_i}{\beta_i} \left(  \sum_{i =1}^N  \frac{1}{\beta_i} -1 \right) (\E \lambda_{t} - \lambda_{t}^{*} ) \\
	 = & \Theta(N^3) \{\Theta (\frac{1}{tN^3}) + \Theta(\frac{1}{(Y^*)^2N^2t^2})\} +
	 \Theta(N^2) \Theta(\frac{1}{N^2t}) \\
	=& \Theta (\frac{1}{t}),
	\end{aligned}
	\end{equation}
	which means that the cumulative regret $R = \sum_t R_t$ is again $\Theta (\log T)$.
\end{proof}

\subsection{Online learning algorithm}

In Section \ref{online}, the online algorithm is shown as a flow chart in Fig. \ref{fig_algo1}. For interested readers, we present the online procedure in its algorithmic form as in Algo. \ref{algo1}. In Algo. \ref{algo1}, the iterated linear regression (Algorithm \ref{lin1}) is repeated at each time instance until the end of the demand response program.

	\begin{algorithm}[h]
		\caption{Online learning algorithm}
		\label{algo1}
		\begin{algorithmic}[1]
			\State \textbf{Initialization}: index of iteration $t=1$, and initial value for the reward signal $\lambda_{0}$.
			\Repeat
			\State In time slot $t$, the operator updates the estimation of parameter $\hat{\gamma}_1^t $ and $\hat{\gamma}_2^t $ using linear regression in \textbf{Algorithm \ref{lin1}}, and
			updates the reward signal
			$$\lambda_{t} = \frac{Y^* d_t + \hat{\gamma}_2^t }{N \hat{\gamma}_{1}^t + N}$$
			and broadcasts to users.
			\State User $i$ responds to the reward signal and reveals its demand reduction as
			$$
			\hat{x}_{i}^{t} = \frac{N \lambda_t - \alpha_i}{\beta_i} + \epsilon_i^t.
			$$
			\State $t=t+1$;
			\Until $t=T$.
			\State \textbf{end}
		\end{algorithmic}
	\end{algorithm}

	\subsection{Ridge regression}
	For the first few updates, there is not enough sample points allowing us to do an efficient linear regression. Therefore, we impose a gaussian prior on the parameters that we want to estimate, and impose a large variance on this prior (a large prior indicates a small regularization parameter $\lambda$ in the penalty, and this parameter does not change the analysis on orders). More specifically, we adopt ridge regression \cite{ridge} to deal with the cold start problem in the linear regression model.

	Adding a prior comes at the cost of a small bias in the estimators to handle with the sample inefficiency in the linear regression model.  For example, for a general cost function there are two regressors int the linear regression models whereas no samples are available in $t_1$ to train the model, and only one sample point is available right after $t_1$. A ridge regression with a small regularization parameter $\lambda$ will validate the least square estimator by assuming a prior information on the estimators. A small $\lambda$ only introduce a negligible bias into the system which avoids the problem of singularity at $t_1$.

	A ridge regression for a linear regression $\mathbf{Z} = \mathbf{X} \bm{\gamma} + \bm{\epsilon}$ is as the follows:
	\begin{equation}
	\hat{\bm{\gamma}} =  (\mathbf{X}^{\text{T}}\mathbf{X} + \lambda I )^{-1}\mathbf{X}^{\text{T}}\mathbf{Z}
	\end{equation}

	As we can see, a small $\lambda$ will avoid the singularity of the matrix $\mathbf{X}^{\text{T}}\mathbf{X}$ and as we gather more sample points, its influence is overwhelmed by $\mathbf{X}^{\text{T}}\mathbf{X}$ and the bias goes down quickly.

\end{appendix}

\bibliographystyle{IEEEtran}	
\bibliography{IEEEabrv,paperref}

\begin{thebibliography}{10}
\providecommand{\url}[1]{#1}
\csname url@samestyle\endcsname
\providecommand{\newblock}{\relax}
\providecommand{\bibinfo}[2]{#2}
\providecommand{\BIBentrySTDinterwordspacing}{\spaceskip=0pt\relax}
\providecommand{\BIBentryALTinterwordstretchfactor}{4}
\providecommand{\BIBentryALTinterwordspacing}{\spaceskip=\fontdimen2\font plus
\BIBentryALTinterwordstretchfactor\fontdimen3\font minus
  \fontdimen4\font\relax}
\providecommand{\BIBforeignlanguage}[2]{{%
\expandafter\ifx\csname l@#1\endcsname\relax
\typeout{** WARNING: IEEEtran.bst: No hyphenation pattern has been}%
\typeout{** loaded for the language `#1'. Using the pattern for}%
\typeout{** the default language instead.}%
\else
\language=\csname l@#1\endcsname
\fi
#2}}
\providecommand{\BIBdecl}{\relax}
\BIBdecl

\bibitem{PalenskyEtAl2011}
P.~Palensky and D.~Dietrich, ``Demand side management: Demand response,
  intelligent energy systems, and smart loads,'' \emph{IEEE transactions on
  industrial informatics}, vol.~7, no.~3, pp. 381--388, 2011.

\bibitem{AlbadiEtAl2008}
M.~H. Albadi and E.~El-Saadany, ``A summary of demand response in electricity
  markets,'' \emph{Electric power systems research}, vol.~78, no.~11, pp.
  1989--1996, 2008.

\bibitem{NguyenEtAl2011}
D.~T. Nguyen, M.~Negnevitsky, and M.~De~Groot, ``Pool-based demand response
  exchange—concept and modeling,'' \emph{IEEE Transactions on Power Systems},
  vol.~26, no.~3, pp. 1677--1685, 2011.

\bibitem{ZhongEtAl2013}
H.~Zhong, L.~Xie, and Q.~Xia, ``Coupon incentive-based demand response: Theory
  and case study,'' \emph{IEEE Transactions on Power Systems}, vol.~28, no.~2,
  pp. 1266--1276, 2013.

\bibitem{ENERNOC2017}
{ENERNOC}, ``Create new revenue with demand response,''
  {\url{https://www.enernoc.com/products/businesses/capabilities/demand-response}},
  2017.

\bibitem{ChenEtAl2014}
C.~Chen, J.~Wang, and S.~Kishore, ``A distributed direct load control approach
  for large-scale residential demand response,'' \emph{IEEE Transactions on
  Power Systems}, vol.~29, no.~5, pp. 2219--2228, 2014.

\bibitem{LiEtAl2011}
N.~Li, L.~Chen, and S.~H. Low, ``Optimal demand response based on utility
  maximization in power networks,'' in \emph{Power and Energy Society General
  Meeting, 2011 IEEE}.\hskip 1em plus 0.5em minus 0.4em\relax IEEE, 2011, pp.
  1--8.

\bibitem{Ohmconnect2017}
{Ohmconnect}, ``Demand response programs,''
  {\url{https://www.ohmconnect.com/wiki/learn-more}}, 2017.

\bibitem{LaiEtAl1985}
T.~L. Lai and H.~Robbins, ``Asymptotically efficient adaptive allocation
  rules,'' \emph{Advances in applied mathematics}, vol.~6, no.~1, pp. 4--22,
  1985.

\bibitem{Qdr2006}
{Department of Energy}, ``Benefits of demand response in electricity markets
  and recommendations for achieving them,'' \emph{Report of the United States
  Congress}, 2006.

\bibitem{SaadEtAl2012}
W.~Saad, Z.~Han, H.~V. Poor, and T.~Basar, ``Game-theoretic methods for the
  smart grid: An overview of microgrid systems, demand-side management, and
  smart grid communications,'' \emph{IEEE Signal Processing Magazine}, vol.~29,
  no.~5, pp. 86--105, 2012.

\bibitem{VardakasEtAl2015}
J.~S. Vardakas, N.~Zorba, and C.~V. Verikoukis, ``A survey on demand response
  programs in smart grids: Pricing methods and optimization algorithms,''
  \emph{IEEE Communications Surveys \& Tutorials}, vol.~17, no.~1, pp.
  152--178, 2015.

\bibitem{TsuiEtAl2012}
K.~M. Tsui and S.-C. Chan, ``Demand response optimization for smart home
  scheduling under real-time pricing,'' \emph{IEEE Transactions on Smart Grid},
  vol.~3, no.~4, pp. 1812--1821, 2012.

\bibitem{JiangEtAl2011}
L.~Jiang and S.~Low, ``Multi-period optimal energy procurement and demand
  response in smart grid with uncertain supply,'' in \emph{Decision and Control
  and European Control Conference (CDC-ECC), 2011 50th IEEE Conference
  on}.\hskip 1em plus 0.5em minus 0.4em\relax IEEE, 2011, pp. 4348--4353.

\bibitem{SamadiEtAl2013}
P.~Samadi, H.~Mohsenian-Rad, V.~W. Wong, and R.~Schober, ``Tackling the load
  uncertainty challenges for energy consumption scheduling in smart grid,''
  \emph{IEEE Transactions on Smart Grid}, vol.~4, no.~2, pp. 1007--1016, 2013.

\bibitem{ZhangEtAl2016}
D.~Zhang, S.~Li, M.~Sun, and Z.~O’Neill, ``An optimal and learning-based
  demand response and home energy management system,'' \emph{IEEE Transactions
  on Smart Grid}, vol.~7, no.~4, pp. 1790--1801, 2016.

\bibitem{ONeillEtAl2010}
D.~O'Neill, M.~Levorato, A.~Goldsmith, and U.~Mitra, ``Residential demand
  response using reinforcement learning,'' in \emph{Smart Grid Communications
  (SmartGridComm), 2010 First IEEE International Conference on}.\hskip 1em plus
  0.5em minus 0.4em\relax IEEE, 2010, pp. 409--414.

\bibitem{XuEtAl2016}
Z.~Xu, T.~Deng, Z.~Hu, Y.~Song, and J.~Wang, ``Data-driven pricing strategy for
  demand-side resource aggregators,'' \emph{IEEE Transactions on Smart Grid},
  2016.

\bibitem{BaltaogluEtAl2016}
S.~Baltaoglu, L.~Tong, and Q.~Zhao, ``Online learning and pricing for demand
  response in smart distribution networks,'' in \emph{Statistical Signal
  Processing Workshop (SSP), 2016 IEEE}.\hskip 1em plus 0.5em minus 0.4em\relax
  IEEE, 2016, pp. 1--5.

\bibitem{KhezeliEtAl2016}
K.~Khezeli and E.~Bitar, ``Risk-sensitive learning and pricing for demand
  response,'' \emph{arXiv preprint arXiv:1611.07098}, 2016.

\bibitem{AlbertEtAl2015a}
A.~Albert and R.~Rajagopal, ``Thermal profiling of residential energy use,''
  \emph{IEEE Transactions on Power Systems}, vol.~30, no.~2, pp. 602--611, Mar.
  2015.

\bibitem{YangEtAl2014}
L.~Yang, D.~Callaway, and C.~Tomlin, ``Dynamic contracts with partial
  observations: application to indirect load control,'' in \emph{American
  Control Conference (ACC)}.\hskip 1em plus 0.5em minus 0.4em\relax IEEE, 2014,
  pp. 1224--1230.

\bibitem{YangEtAl2015}
------, ``Indirect load control for electricity market risk management via
  risk-limiting dynamic contracts,'' in \emph{American Control Conference
  (ACC)}.\hskip 1em plus 0.5em minus 0.4em\relax IEEE, 2015, pp. 3025--3031.

\bibitem{Shaley2012}
S.~Shalev-Shwartz, ``Online learning and online convex optimization,''
  \emph{Foundations and Trends in Machine Learning}, vol.~4, no.~2, pp.
  107--194, 2012.

\bibitem{Bubeck2012}
S.~Bubeck and N.~Cesa-Bianchi, ``Regret analysis of stochastic and
  nonstochastic multi-armed bandit problems,'' \emph{Foundations and Trends in
  Machine Learning}, vol.~5, no.~1, pp. 1--122, 2012.

\bibitem{Chow76}
G.~P. Chow, \emph{Analysis and Control of Dynamic Economic Systems}.\hskip 1em
  plus 0.5em minus 0.4em\relax New York: Wiley, 1976.

\bibitem{boyd2004convex}
S.~Boyd and L.~Vandenberghe, \emph{Convex optimization}.\hskip 1em plus 0.5em
  minus 0.4em\relax Cambridge university press, 2004.

\bibitem{BLUE}
C.~Henderson, ``Best linear unbiased estimation and prediction under a
  selection model,'' \emph{Biometrics}, pp. 423--447, 1975.

\bibitem{algo}
T.~Cormen, \emph{Introduction to algorithms}.\hskip 1em plus 0.5em minus
  0.4em\relax MIT press, 2009.

\bibitem{Keskin2014}
N.~Keskin and A.~Zeevi, ``Dynamic pricing with an unknown demand model:
  Asymptotically optimal semi-myopic policies,'' \emph{Operations Research},
  vol.~62, no.~5, pp. 1142--1167, 2014.

\bibitem{ridge}
N.~Draper, H.~Smith, and E.~Pownell, \emph{Applied regression analysis}.\hskip
  1em plus 0.5em minus 0.4em\relax New York: Wiley, 1966, vol.~3.

\bibitem{ratio}
H.~Seltman, ``Approximations for mean and variance of a ratio,'' 2012,
  unpublished note.

\end{thebibliography}

\end{document}